\documentclass[11pt]{article}
\usepackage{amssymb}
\usepackage{latexsym}
\usepackage{epsfig}

\usepackage{amstex}
\usepackage{amsxtra}
\usepackage{amstext}
\usepackage{amsfonts}
\usepackage{amsbsy}
\usepackage{amsopn}
\usepackage{amssymb}

\begin{document}
\newtheorem{definition}{Definition}
\newtheorem{theorem}{Theorem}
\newtheorem{fact}{Fact}
\newtheorem{lemma}{Lemma}
\newtheorem{proposition}{Proposition}
\newtheorem{remark}{Remark}
\newtheorem{corollary}{Corollary}
\newtheorem{conjecture}{Conjecture}
\newtheorem{example}{Example}
\input amssym.def
\input amssym
\def\nin{\noindent}

\title{{\bf Variation on a theme of Selberg integrals}}
%Subject Classification.} Primary 58F05. Secondary }%}
\author{ A. Kazarnovski-Krol} 
Department of Mathematics\\  Rutgers University \\
New Brunswick , NJ 08854, USA\\
        {\sf email: akrol \@math.rutgers.edu}
%\date{August, 1996}
\maketitle

\begin{abstract} In this paper we calculate some Generalized Selberg
integrals. The answer is expressed in terms of
$\Gamma$-functions. Integrals of this type serve as normalization constants or
directly via undoing  2-D integrals for determination of structural
constants of operator algebra. 
\end{abstract}
\section{Notations}
$n, N$ two positive integers

 $\alpha_1, \alpha_2, \ldots, \alpha_{n-1}$ simple roots of root
system of type $A_{n-1}$

 $\alpha_1, \alpha_2, \ldots , \alpha_{N-1}$ simple roots of root
system of type $A_{N-1}$

$\Sigma_{+}(N-1)$ positive roots of root system of type $A_{N-1}$

$$\rho = \frac{1}{\kappa} \sum_{\alpha \in \Sigma_{+}(N-1)} \alpha$$

$$q=e^{\frac{2 \pi i}{\kappa}}$$

\section{Introduction} 
The paper is a continuation of the series \cite{A1}, \cite{A2},
\cite{A3}, \cite{A4}, \cite{A5}  
and contains an extension of the results of \cite{A3}.

While in the previous papers we used conformal field theory for the
needs of harmonic analysis ( $(n+1)!$ dimensional space of solutions
to the Heckman-Opdam hypergeometric system (wave functions in
Calogero-Sutherland model) related to root system of type $A_n$
is isomorphic to certain $(n+1)!$ dimensional space of conformal blocks
in $WA_n$ algebra, Harish-Chandra asymptotic solutions provide a basis
in the space of conformal blocks, zonal spherical function is a
particular conformal block) , here the accent is changed to
applications to conformal field theory.
The objective of conformal field theory is to construct Green's
function ,i.e. monodromy invariant function which depends on
$z=(z_1,z_2,\ldots ,z_N)$ and conjugate variables.
Green's function is constructed out of conformal blocks:
$$G = \sum c_i |F_i(z)|^2$$
Conformal blocks are asymptotic solutions with the prescribed
asymptotic behaviour at the singularity.
 Asymptotic solutions are normalized so that 
the leading asymptotic coefficient is equal to $1$.
On the other hand asymptotic solutions (conformal blocks) are provided
by certain multidimensional integrals, and in order to calculate leading
asymptotic coefficients value of Selberg-type integrals is needed.
Also the Green's function can be obtained by undoing $2 D$
 integrals cf. \cite{Do}. 
The coefficients $c_i$ should be properly normalized s.t. $<\phi, \phi>=1$.
 In this paper we calculate Selberg-type integrals related to the
$$V_{\lambda} \otimes V_{N \Lambda_1} \longrightarrow
V_{\lambda^{\prime}}.$$
Parameter $\lambda$ is assumed to be generic, $\Lambda_1$ is the first
fundamental weight. 

Note that in this case there is no multiplicities in the tensor
product and the answer is expressed as a product of $\Gamma$-functions. 

Here is organization of the paper.
In theorem \ref{15t1} we prove that certain type of integrals
satisfy Heckman-Opdam hypergeometric system. 

Then we collapse the arguments to the unity and using Opdam's result
obtain in theoem  \ref{15Sel} the value of Selberg-type integrals .
Here is organization of the paper.
In section \ref{15S1} we introduce the necessary combinatorics including two
root systems. The first one  is related to variables of integration,
while the second one is related to the second factor of tensor product
 $N \Lambda_1$ (the row with $N$ boxes).

Then in section 4 we  compare results  with the usual Selberg integral
, and finally in section 5 we emphasize one particular case , when after
integration we actually obtain a monomial.
 
\section{}
\label{15S1}
Fix two positive integers $n$ and $N$. We distinguish two root systems
$A_{n-1}$ and $A_{N-1}$ with simple roots $\alpha_1, \alpha_2, \ldots,
\alpha_{n-1}$ and $\alpha_1, \alpha_2, \ldots , \alpha_{N-1}$,
correspondigly.  It is not important whether $n < N$ or $n \ge N$.
We would like to realize this situation as follows. Let $m$ be the
maximum of $N$ and $n$ : 
$$m={\rm Maximum}(n,N).$$ Consider $m$-dimensional Euclidean vector space  
$\mathbb R^m$ with $e_1, e_2, \ldots , e_m$ as orthonormal basis.
Let $\alpha_1=e_1-e_2, \alpha_2=e_2-e_3, \ldots,
\alpha_{m-1}=e_{m-1}-e_m$. The set of positive roots of $A_{N-1}$ shall
be denoted by $\Sigma_{+}(N-1)$.

Let $\Lambda_1, \Lambda_2 , \ldots, \Lambda_{m-1}$ be fundamental
weights, i.e. $(\Lambda_i, \alpha_j)=\delta_{ij}$ and $\Lambda_i \in
span\{\alpha_j\}$, i.e.
$(\Lambda_i, e_1 +e_2 +\ldots +e_m)=0$.

Let $z=(z_1, z_2, \ldots, z_N)$. We assume that $z_1,z_2, \ldots, z_N$
are real and $$0 < z_1 < z_2 < \ldots < z_N.$$
Let also $h_1=\Lambda_1, h_2=\Lambda_1-\alpha_1, \ldots ,
h_{n+1}=\Lambda_1-\alpha_1-\ldots -\alpha_{n}$ .

Let $V_{\lambda}$ denotes the irreducible highest weight module  over
quantum group $U_q(sl(n+1))$ with highest weight vector $\lambda$.
$$q=e^{\frac{2 \pi i}{\kappa}}$$ 

Consider
\begin{equation}
V_{\lambda}\stackrel{\otimes V_{\Lambda_1}}\longrightarrow
V_{\lambda+h_{i_1}}
\stackrel{\otimes V_{\Lambda_1}}\longrightarrow
V_{\lambda+h_{i_1} +h_{i_2}}\stackrel{\otimes V_{\Lambda_1}}
\longrightarrow \ldots \stackrel{\otimes V_{\Lambda_1}}\longrightarrow
V_{\lambda+h_{i_1}+ h_{i_2} +\ldots+ h_{i_N}}
\label{13Veq0}
\end{equation}

Essentially we fix the ordered set of $N$ numbers ( each of them not
exceeding $n$) .

\begin{equation}
I=\left\{ i_1, i_2, \ldots , i_N \right \}.
\label{}
\end{equation}

Let 
\begin{equation}
s(k)=\# \Big\{ i_j \in I |\; i_j > k \Big\},
\label{}
\end{equation}
where by $\#$  we denote the cardinality of the corresponding set.
Without loss of generality we assume 
\begin{equation}
n={\rm Maximum}\Big\{ i_j-1 | \; i_j \in I \; \Big \}.
\label{}
\end{equation}
Let 
\begin{equation}
a(p)=\#\Big\{ i_j >1 | j > p \Big\}.
\label{a}
\end{equation}

Consider then the set of variables 
$$\left\{t_{ij} | i=1, \ldots, s(j), j= 1 , \ldots,n \right\}$$
To each variable $t_{ij}$ we assign a simple root $\alpha_j$.
To variables $z_i$ we assign the first fundamental weight $\Lambda_1$,
also we fix $z_0=0$ and assign to $z_0$  vector $\lambda$, $\lambda
\in span\{\alpha_i |i=1, \ldots, n-1\}$.

Accordingly, we consider the integral
\begin{eqnarray}
\int \prod t_{ij}^{\frac{(\lambda, -\alpha_j)}{\kappa}}
\prod (z_i -t_{i1})^{\frac{(\Lambda_1, -\alpha_1)}{\kappa}}
\prod (t_{i_1 j_1}-t_{i_2 j_2})^{\frac{(\alpha_{j_1},\alpha_{j_2})}
{\kappa} }\frac{dt_{11}}{t_{11}} \ldots \frac{dt_{s(n),n}}{t_{s(n),n}}
\label{13Veq1} 
\end{eqnarray}

\begin{proposition}
\label{15PR1}

The leading asymptotic of the integral (\ref{13Veq1}) over the cycle
related to conformal block (\ref{13Veq0}) is equal:
\begin{eqnarray}
z^{\mu}=z_1^{\mu_1} z_2^{\mu_2} \ldots z_N^{\mu_N}\nonumber \\
=
\prod z_j^{\frac{(\lambda+h_{i_1}+h_{i_2} + \ldots + h_{i_{j-1}},
\sum_{i=1}^{i_j -1} - \alpha_i)-(i_j-1) +(h_{i_1} + h_{i_2} + \ldots +
h_{i_{j-1}} -(j-1)\Lambda_1, \Lambda_1) }{\kappa}}
\label{13Veq2}
\end{eqnarray}
\begin{equation}
\mu_j=\frac{(\lambda+h_{i_1} +\ldots +h_{i_{j-1}}, \sum_{i=1}^{i_j-1} 
 - \alpha_i)-(i_j-1) +(h_{i_1}+h_{i_2}+ \ldots + h_{i_{j-1}} -(j-1)
\Lambda_1,\Lambda_1)}{\kappa} \label{13Veq2.2}
\end{equation}
The leading asymptotic coefficient is equal:

\begin{eqnarray}
\prod_{j=1}^N \prod_{p=1}^{i_j -1}
exp( \pi i (\frac{(\lambda+h_{i_1}+h_{i_2} +\ldots+ h_{i_{j-1}},
\sum_{i=i_j -p}^{i_j-1} -\alpha_i)}{\kappa} -\frac{p-1}{\kappa}
))\nonumber \\
\times \frac{\Gamma(1-\frac{1}{\kappa}) (2 \pi i)}
{\Gamma(\frac{(\lambda +h_{i_1} +\ldots +
h_{i_{j-1}},\sum_{i=i_j-p}^{i_j-1} -\alpha_i)}{\kappa}  -
\frac{p}{\kappa} +1) 
\Gamma(\frac{(\lambda + h_{i_1}+ h_{i_2} + \ldots + h_{i_j -1},
\sum_{i_j-p}^{i_j-1} \alpha_i)}{\kappa} + \frac{p-1}{\kappa}+1)}
\label{13Veq3}
\end{eqnarray}
Moreover , for generic $\lambda, \kappa$
\begin{equation}
\mu_i -\mu_j+ \frac{1}{\kappa}(j-i) \notin \mathbb Z
\end{equation}

\end{proposition}
\begin{definition}
Define $\rho$ to be $\frac{1}{\kappa}$ times half the sum of positive
roots of root system of type $A_{N-1}$ :

\begin{eqnarray}
\rho =(\rho_1, \rho_2, \ldots, \rho_N)\nonumber \\
= \frac{1}{2 \kappa}(N-1, N-3, \ldots, 3-N, 1-N)
\label{13Veq4}
\end{eqnarray}
\end{definition}
{\bf Proof:}
 The calculation of the leading asymptotic coefficient inductively
uses Dirichlet's formula (and essentially corresponds to composition of
vertex operators).
Let 
\begin{eqnarray}
\eta=\mu-\rho \\
\eta_j=\mu_j -\rho_j
\label{13Veq5}
\end{eqnarray}

We want check that:
\begin{equation}
\eta_s-\eta_j \notin \mathbb Z
\label{15V1}
\end{equation}
 for $s \ne j$.

We consider the two cases:1) $i_s \ne i_j$ and 2) $i_s=i_j$.

1). In the first case : $i_s \ne i_j $ we get
\begin{equation}
\eta_s -\eta_j =
\frac{\lambda_{i_s} -\lambda_{i_j}}{\kappa} +\ldots 
\label{15Veq2}
\end{equation}
Thus (\ref{15V1}) is satisfied for generic $\lambda$.

2). In the second case let $i_s=i_j=p+1$
\begin{eqnarray}
\eta_j -\eta_s=
\frac{1}{\kappa}\Big((h_{p+1} +h_{i_{s+1}} +h_{i_{s+2}} + \ldots
+h_{i_{j-1}}, -(e_1-e_{p+1}))\nonumber \\+
(h_{p+1} + h_{i_{s+1}} + \ldots + h_{i_{j-1}} -(j-s) \Lambda_1,
\Lambda_1) +(j-s) \Big) \ne 0
\label{15Veq3}
\end{eqnarray}
As $(h_{p+1}, -\alpha_1-\ldots -\alpha_{p})=1$.
Moreover, (\ref{15Veq3}) is positive for positive $\kappa$ and thus the
singularities of (\ref{15V1}) can be avoided by assuming $\kappa$ 
to be irrational.
So from now on, we  work in the above hypotheses  on $\kappa$ and $\lambda$. 
          $$\Box$$
\vskip 1cm

Let $L$ be the following differential operator
\begin{equation}
L=
{\sum_{i=1}^{N} (z_i \frac{\partial} {\partial z_i})^2} -
{k\sum_{i < j} {\frac{z_j +z_i}{z_j-z_i}}
\big(z_i \frac{\partial} {\partial z_i}-
z_j \frac{\partial} {\partial z_j}\big)}.
\end{equation}

Explicitly, in the case of root system of type $A_{N-1}$ the generators of
the system of commuting operators were calculated in ref \cite{Sek}.
Namely, let
\begin{eqnarray}
D(\zeta,k)= \frac{1}{\prod_{i<j}(z_i-z_j)}\sum_{w \in S_{N}}
det(w) \prod_{j=1}^{N} z_j^{N-w(j)} \prod_{i=1}^{N}
(\zeta + z_i \frac{\partial}{\partial z_i}+(w \delta, e_i) k)
\end{eqnarray}
Then 
\begin{equation}
D(\zeta, k)=\sum_{r=1}^{N} \zeta^r D^{(k)}_{N-r}
\end{equation}
Operators $D^{(k)}_{N-r}$  commute with each other.
Moreover, the system of hypergeometric differential equations can be
written in the form:

\begin{equation}
D(\zeta,k) \phi(\eta,k,z)=\prod_{i=1}^{N}(\zeta +\eta_i)
\phi(\eta,k,z)
\label{syst}
\end{equation}

Consider  solutions of
\begin{equation}
L \phi =((\eta , \eta) -(\rho, \rho)) \phi
\end{equation}
 
of the form

\begin{equation}
\phi= \phi(\mu,k,z)=z^{\mu}(1 +\ldots)
\quad ,
\end{equation}

where  $\mu = w \eta +\rho$.

Then for
generic $\eta$ $((\eta, \alpha^{\vee}) \notin \Bbb Z, \;
\alpha \in \Sigma)$ solutions $\phi(w \eta +\rho,k,z)$ are
uniquely defined, converge to analytic function in chamber $0 < |z_1| <
|z_2| < \ldots < |z_{N}|$,
linearly independent and provide a basis of linear space of solutions
to the whole system of hypergeometric equations (\ref{syst}).

These solutions $\phi(w \eta+\rho,k,z)$ will be referred to as {\bf{Harish
Chandra asymptotic solutions}}.

\begin{theorem}(Opdam)
The value of asymptotic solution $\phi(\eta,k,z)$ at the unity is equal:

\begin{eqnarray}
\phi(\eta,k,1)=\lim_{z \to 1} \phi(\eta,k,z)=\nonumber \\
\frac{\prod_{\alpha \in \Sigma_{+}(N-1)}
\frac{\Gamma((\eta, \alpha)+1)}
     {\Gamma((\eta, \alpha)- \frac{1}{\kappa}+1)}  }
{\prod_{\alpha \in \Sigma_{+}(N-1)}
\frac{\Gamma(-(\rho,\alpha)+1)}
     {\Gamma(-(\rho,\alpha)-\frac{1}{\kappa}+1)}} 
\label{13Veq6}
\end{eqnarray}
cf. \cite{Op1}, theorem 6.3.
\label{15tt}
\end{theorem}

\begin{remark}
Opdam's result is actually in more general context, see also
\cite{A1},\cite{A2}, \cite{A4}, \cite{A5}.
\end{remark}

\begin{remark}
The structure of the constant in theorem \ref{15tt} is still very much
similar to those of  $c$-function of Harish-Chandra. It is not
accidental, since it is obtained essentially using Harish Chandra
decomposition for zonal spherical function and monodromy properties of
zonal spherical function. Recall that $c$-function of Harish Chandra
is needed for the Plancherel measure
$$d \mu = \frac{d \lambda}{|c(\lambda,k)|^2}.$$
$c$-function was introduced by Harish Chandra, in the case of $SL(n,\mathbb C)$
it was calculated by Gelfand and Naimark, in the case of $SL(n,\mathbb R)$
by Bhanu Murti, and in general case by Gindikin and Karpelevich.
\end{remark}

\begin{theorem}
The following integrals

\begin{eqnarray}
\int \prod t_{ij}^{\frac{(\lambda, -\alpha_j)}{\kappa}}
\prod (z_i -t_{i1})^{\frac{(\Lambda_1, -\alpha_1)}{\kappa}}
\prod (t_{i_1 j_1}-t_{i_2 j_2})^{\frac{(\alpha_{j_1},\alpha_{j_2})}
{\kappa} }\frac{dt_{11}}{t_{11}} \ldots \frac{dt_{s(n),n}}{t_{s(n),n}}
\label{13Veqq} 
\end{eqnarray}
over appropriate cycles is a common eigenfunction of hypergeometric
system of differential equations (\ref{syst}).
Cycles should be chosen so that they are encoded  as  singular vectors of
the tensor product of irreducible highest weight 
modules over quantum group $U_q(sl(n+1))$ :
$$ 
V_{\lambda} \otimes V_{\Lambda_1} \otimes V_{\Lambda_1} \otimes \ldots
\otimes V_{\Lambda_1}
$$
\label{15t1} 
\end{theorem}
{\bf Proof:}
 The theorem is proved with the help of integration by parts
``upwards''. It is convenient to write  the integral in the
hierarchical form, and to organize variables of integration  in the manner
of Gelfand-Tsetlin patterns, s.t. the number of variables in the upper
row is greater than or equal to the number of variables in the
lower row. Since $(\alpha_i, \alpha_j)=0$ if $|i-j| >1$ and
$(\Lambda_1,\alpha_i)= \delta_{1i}$,
there is only interaction between the variables in the adjacent rows. 
It is quite instructive to think over  what happens if the number of
variables in the next row is the same as in the previous one.
 In this respect see example \ref{15exam1}  below. 
In the same way we  obtained an integral representation   in \cite{A4}.
                                                                $\Box$
\begin{remark}
Integrals of the type considered in the theorem \ref{15t1}
are produced as conformal blocks of $W_n$ algebra, except that we do
not have the factor $\prod_{i < j} (z_i -
z_j)^{\frac{(\Lambda(i),\Lambda(j))}
{\kappa}}$
before the integral. One should notice, that this factor does not
change the homological structure of the fiber , but affects
Gauss-Manin connection.
\end{remark}
\begin{remark}
Passage from ''loops'' to ``tines'' kills the kernel of the  contravariant
form and adds quantum Serre's relations \cite{SV1},\cite{Va1}.
\end{remark}

\begin{example} 
\label{15exam1}
Triple integral for usual hypergeometric function.
We consider the following situation:
$$
V_{\lambda} \stackrel{\otimes V_{\Lambda_1}} \longrightarrow
V_{\lambda+h_3} \stackrel{\otimes V_{\Lambda_1}} \longrightarrow
V_{\lambda+h_3+h_2}
$$

\begin{eqnarray}
\int \int \int
(t_1 t_2)^{\frac{(\lambda, -\alpha_1)}{\kappa}}
(z_1-t_1)^{-\frac{1}{\kappa}}(z_2-t_1)^{-\frac{1}{\kappa}}
(z_2-t_1)^{-\frac{1}{\kappa}} (z_2-t_2)^{-\frac{1}{\kappa}}
(t_2-t_1)^{\frac{2}{\kappa}} \nonumber \\
\times
t_3^{\frac{(\lambda,-\alpha_2)}{\kappa}} (t_1-t_3)^{-\frac{1}{\kappa}}
(t_2-t_3)^{-\frac{1}{\kappa}} \frac{d t_1}{t_1} \frac{d t_2}{t_2}
\frac{dt_3}{t_3} \nonumber \\=
const  \; z_1^{\frac{(\lambda,-\alpha_1-\alpha_2)}{\kappa} -\frac{2}{\kappa}}
z_2^{\frac{(\lambda,-\alpha_1)}{\kappa} -\frac{2}{\kappa}}
F(\frac{1}{\kappa}, \frac{(\lambda, -\alpha_2)}{\kappa},
\frac{(\lambda, -\alpha_2)}{\kappa} -\frac{1}{\kappa} +1,
\frac{z_1}{z_2})
\label{15pp1}
\end{eqnarray}

As a corollary we obtain the following identity involving 4-indices
summation:
\begin{eqnarray}
\frac{\Gamma(\frac{(\lambda, \alpha_1)}{\kappa}+1)}
{\Gamma( \frac{(\lambda, \alpha_1)}{\kappa} +\frac{1}{\kappa})}
\frac{\Gamma(\frac{(\lambda, -\alpha_2)}{\kappa} -\frac{1}{\kappa}+1)}
{\Gamma(\frac{(\lambda, -\alpha_2)}{\kappa})}
\frac{\Gamma(\frac{(\lambda, -\alpha_1
-\alpha_2)}{\kappa}+1-\frac{2}{\kappa})}
{\Gamma(\frac{(\lambda, -\alpha_1 -\alpha_2)}{\kappa}
-\frac{1}{\kappa})}
\frac{1}
{\Gamma(\frac{1}{\kappa}) \Gamma(\frac{1}{\kappa})
\Gamma(\frac{-2}{\kappa})
\Gamma(\frac{1}{\kappa})}
 \nonumber \\ \times
\sum_{m_1+m_2+m_3+m_4=M} \quad
\frac{\Gamma(\frac{(\lambda, -\alpha_1-\alpha_2)}{\kappa} +m_1 +m_3
+m_4 -\frac{1}{\kappa})}  
{\Gamma(\frac{(\lambda,-\alpha_1-\alpha_2)}{\kappa} +m_1 +m_3 +m_4 +1-
\frac{2}{\kappa})}  \nonumber \\
\times
\frac{\Gamma(\frac{(\lambda, -\alpha_2)}{\kappa} +m_4)}
{\Gamma(\frac{(\lambda,-\alpha_2)}{\kappa}+m_4+1-\frac{1}{\kappa})}\nonumber
\\ \times 
\frac{\Gamma(\frac{(\lambda, \alpha_1)}{\kappa} + \frac{1}{\kappa}+m_2
+m_3 +m_4)}
{\Gamma(\frac{(\lambda, \alpha_1)}{\kappa} +1 +m_2 +m_3
+m_4)}\nonumber \\ 
\times
\frac{
\Gamma(\frac{1}{\kappa}+m_1) \Gamma(\frac{1}{\kappa}+m_2)
\Gamma(\frac{-2}{\kappa}+m_3) \Gamma(\frac{1}{\kappa}+m_4)}
{ m_1 ! m_2 ! m_3 ! m_4 !}
\nonumber \\
=
\frac{\Gamma(\frac{(\lambda, -\alpha_2)}{\kappa} -\frac{1}{\kappa}+1)}
{\Gamma(\frac{1}{\kappa}) \Gamma(\frac{(\lambda, -\alpha_2)}{\kappa})} \times
\frac{\Gamma(\frac{1}{\kappa}+M) \Gamma(\frac{(\lambda,
-\alpha_2)}{\kappa}+M)}
{\Gamma(\frac{(\lambda, -\alpha_2)}{\kappa} -\frac{1}{\kappa}+1+M) M !}
\label{}
\end{eqnarray}
\end{example}

Here is the main result of the paper.

\begin{theorem} The Generalized Selberg integral is equal:
\label{15Sel}
\begin{eqnarray}
\int
\prod {t_{ij}}^{\frac{(\lambda, -\alpha_j)}{\kappa}}   
\prod (1-t_{i1})^{-\frac{1}{\kappa}}
\prod (t_{i_1 j_1 }-t_{i_2 j_2})^{\frac{(\alpha_{j_1},\alpha_{j_2})}{\kappa}}
\frac{d t_{11}}{t_{11}} \ldots \frac{d t_{s(n),n}}{t_{s(n),n}}\nonumber \\  
=
A \prod_{j=1}^N \prod_{p=1}^{i_j -1}
exp( \pi i (\frac{(\lambda+h_{i_1}+h_{i_2} +\ldots + h_{i_{j-1}},
\sum_{i=i_j -p}^{i_j-1} -\alpha_i)}{\kappa} -\frac{p-1}{\kappa}
))\nonumber \\
\times \frac{\Gamma(1-\frac{1}{\kappa}) (2 \pi i)}
{\Gamma(\frac{(\lambda +h_{i_1} + \ldots +
h_{i_{j-1}},\sum_{i=i_j-p}^{i_j-1} -\alpha_i)}{\kappa}  -
\frac{p}{\kappa} +1) 
\Gamma(\frac{(\lambda+h_{i_1}+h_{i_2} +\ldots + h_{i_j -1},
\sum_{i_j-p}^{i_j-1} \alpha_i)}{\kappa} + \frac{p-1}{\kappa}+1)}
\nonumber \\
\times
\frac{\prod_{\alpha \in \Sigma_{+}(N-1)}
\frac{ \Gamma((\mu - \rho, \alpha)+1)}
     { \Gamma((\mu-\rho, \alpha)- \frac{1}{\kappa}+1)} }
{\prod_{\alpha \in \Sigma_{+}(N-1)}
\frac{\Gamma(-(\rho,\alpha)+1)}
     {\Gamma(-(\rho,\alpha)-\frac{1}{\kappa}+1)}}
\label{13Veq}
\end{eqnarray}

\begin{equation}
A=e^{-\frac{\pi i}{\kappa} \sum_{p=1}^N a(p)}
\label{A}
\end{equation} 

Here $a(p)$ is defined by (\ref{a}), $\mu$ by formula (\ref{13Veq2.2}), and
$\rho$ by (\ref{13Veq4}).

$$
V_{\lambda} \otimes V_{N \Lambda_1} \longrightarrow V_{\lambda +h_{i_1}
+ h_{i_2} +\ldots + h_{i_N}}
$$

\end{theorem}
{\bf Proof:}
The theorem immediately follows from theorems \ref{15tt}, \ref{15t1}.
The constant $A$ takes into account the phase which is earned as $z_i$ goes through
some $t_{1j}$ as all $z_i$ collapse to the unity.  
$|z_i|< |t_{1j}|$ changes to $|t_{1j}| < |z_i|$.

\begin{example}
Consider $V_{\lambda} \otimes V_{3 \Lambda_1} \longrightarrow
V_{\lambda+2 h_1 +h_2}$.
Accordingly we consider the integral
\begin{equation}
\int t^{\frac{(\lambda, -\alpha)}{\kappa}-1} (1-t)^{-\frac{3}{\kappa}}
dt= \frac{ \Gamma(1-\frac{3}{\kappa}) (2  \pi i) e^{\pi i \frac{(\lambda, -\alpha)}{\kappa}}}
{\Gamma(\frac{\lambda, \alpha}{\kappa}+1) \Gamma( \frac{(\lambda,
-\alpha)}{\kappa}+1 -\frac{3}{\kappa})}
\label{15ex1}
\end{equation}
Here the contour of integration is chosen so that it starts and ends
at $1$ and encloses $0$ anticlockwise.
Now consider the following three cases.

a).
$$
V_{\lambda} \stackrel{ \otimes V_{\Lambda_1}} \longrightarrow
V_{\lambda +h_2} \stackrel{\otimes V_{\Lambda_1}} \longrightarrow
V_{\lambda+h_2 +h_1} \stackrel{\otimes V_{\Lambda_1}} \longrightarrow
V_{\lambda+h_2 +2 h_1}
$$

Then 
$$A=1$$
$$
\mu=( \frac{(\lambda, -\alpha)}{\kappa} -\frac{1}{\kappa},
-\frac{1}{\kappa}, -\frac{1}{\kappa})
$$

$$\rho=(\frac{1}{\kappa}, 0, -\frac{1}{\kappa})$$
$$\mu-\rho=(\frac{(\lambda,-\alpha)}{\kappa} - \frac{2}{\kappa}, 
-\frac{1}{\kappa}, 0)$$
So according to the theorem the Selberg integral in this case is
equal:
\begin{eqnarray}
1 \times 
\frac{e^{\pi i \frac{(\lambda, -\alpha)}{\kappa}}
\Gamma(1-\frac{1}{\kappa})  (2 \pi i)}
{\Gamma(\frac{(\lambda, -\alpha)}{\kappa} +1 -\frac{1}{\kappa})
\Gamma(\frac{(\lambda, \alpha)}{\kappa}+1)} \nonumber \\
\times 
\frac{\Gamma(\frac{(\lambda, -\alpha)}{\kappa} - \frac{2}{\kappa} +1)
\Gamma(\frac{(\lambda, -\alpha)}{\kappa}-\frac{1}{\kappa}+1)
\Gamma(1-\frac{1}{\kappa})}
{\Gamma(\frac{(\lambda, -\alpha)}{\kappa} -\frac{3}{\kappa} +1)
\Gamma(\frac{(\lambda, -\alpha)}{\kappa} - \frac{2}{\kappa}+1)
\Gamma(1-\frac{2}{\kappa})} \over
\frac{\Gamma(1-\frac{1}{\kappa}) \Gamma(1-\frac{1}{\kappa})
\Gamma(1-\frac{2}{\kappa})}
{\Gamma(1-\frac{2}{\kappa}) \Gamma(1-\frac{2}{\kappa}) 
\Gamma(1-\frac{3}{\kappa})}
\label{15ex2}
\end{eqnarray}
After simplifications we get (\ref{15ex1}).

b).
Consider now the following choice of intermediate channels:
$$
V_{\lambda} \stackrel{\otimes V_{\Lambda_1}} \longrightarrow
V_{\lambda+ h_1} \stackrel{ \otimes V_{\Lambda_1}} \longrightarrow
V_{\lambda+h_2 +h_1} \stackrel{\otimes V_{\Lambda_1}} \longrightarrow
V_{\lambda+h_2 +2 h_1}
$$

$$ A=e^{-\frac{\pi i}{\kappa}} $$
$$\mu =(0, \frac{(\lambda, -\alpha)}{\kappa} -\frac{1}{\kappa}, 
-\frac{1}{\kappa})$$
$$\rho =(\frac{1}{\kappa}, 0, -\frac{1}{\kappa})$$

$$\mu -\rho=
(-\frac{1}{\kappa}, \frac{(\lambda, -\alpha)}{\kappa}
-\frac{1}{\kappa}, -\frac{1}{\kappa})
$$

So according to the theorem we get that the value of Selberg integral
is equal:
\begin{eqnarray}
e^{-\frac{\pi i}{\kappa} } \times 
\frac{e^{\pi i \frac{(\lambda+h_1, -\alpha)}{\kappa} }
\Gamma(1-\frac{1}{\kappa}) ( 2 \pi i )}
{\Gamma(\frac{(\lambda, -\alpha)}{\kappa} +1) \Gamma(\frac{(\lambda,
\alpha)}{\kappa}+1)}
\nonumber \\
\times 
\frac{
\frac{\Gamma(\frac{(\lambda, -\alpha)}{\kappa}+1)
\Gamma(\frac{(\lambda, -\alpha)}{\kappa} -\frac{1}{\kappa}+1)
\Gamma(1-\frac{1}{\kappa})}
{\Gamma(\frac{(\lambda, -\alpha)}{\kappa}+1-\frac{2}{\kappa})
\Gamma(1-\frac{2}{\kappa})}}
{
\frac{\Gamma(1-\frac{1}{\kappa}) \Gamma(1-\frac{1}{\kappa}) 
\Gamma(1-\frac{2}{\kappa})}
{\Gamma(1-\frac{2}{\kappa}) \Gamma(1-\frac{2}{\kappa}) 
\Gamma(1-\frac{3}{\kappa})}}
\label{15ex3}
\end{eqnarray}
 Again in agreement with (\ref{15ex1})

c).
$$
v_{\lambda} \stackrel{\otimes V_{\Lambda_1}} \longrightarrow
V_{\lambda +h_1} \stackrel{\otimes V_{\Lambda_1}} \longrightarrow
V_{\lambda + 2 h_1} \stackrel{\otimes V_{\Lambda_1}} \longrightarrow
V_{\lambda+2 h_1 +h_2}
$$

$$A=e^{-\frac{2 \pi i }{\kappa}}$$

$$\mu=(0,0, \frac{(\lambda, -\alpha)}{\kappa} - \frac{3}{\kappa})$$

$$\mu -\rho=(-\frac{1}{\kappa}, 0, \frac{(\lambda, -\alpha)}{\kappa} -\frac{2}{\kappa})$$
So according to the theorem the value of the Selberg integral is
equal:
\begin{eqnarray}
e^{-\frac{2 \pi i }{\kappa}} 
\times
\frac{e^{\pi i \frac{(\lambda+ 2 \Lambda_1, -\alpha)}{\kappa}}
\Gamma(1-\frac{1}{\kappa}) ( 2 \pi i )}
{\Gamma(\frac{(\lambda, -\alpha)}{\kappa} - \frac{3}{\kappa}+1)
\Gamma(\frac{(\lambda, \alpha)}{\kappa} + \frac{2}{\kappa}+1)
}
\nonumber \\
\times
\frac{ 
\frac{
\Gamma(1-\frac{1}{\kappa}) \Gamma(\frac{(\lambda, \alpha)}{\kappa} + 
\frac{2}{\kappa}+1) \Gamma(\frac{(\lambda, \alpha)}{\kappa}
+\frac{1}{\kappa}+1)}
{
\Gamma(1-\frac{2}{\kappa}) \Gamma( \frac{(\lambda, \alpha)}{\kappa} +
\frac{1}{\kappa} +1) \Gamma(\frac{(\lambda,\alpha)}{\kappa}+1)
}
}
{\frac
      {\Gamma(1-\frac{1}{\kappa}) \Gamma(1-\frac{1}{\kappa}) 
       \Gamma(1-\frac{2}{\kappa})}
     {\Gamma(1 -\frac{2}{\kappa}) \Gamma(1-\frac{2}{\kappa})
\Gamma(1-\frac{3}{\kappa})}}
\label{15ex4}
\end{eqnarray}
\end{example}
Again after simplifications we get  (\ref{15ex1})
\begin{remark}
The Selberg type integrals are needed as normalization constants or
directly  undoing the $2D$ integrals for calculation of structural
constants of operator algebra \cite{Do} ,\cite{DF}.
\end{remark}

\section{Comparison with usual Selberg integral}
In ref \cite{Sel} Selberg considered the following multiple integral
\begin{eqnarray} 
\int \ldots \int 
\prod_i t_i^{a-1} (1-t_i)^{b-1} \prod_{i<j}|t_i- t_j|^{2c} dt_1 \ldots
dt_m \nonumber \\=
\prod_{j=0}^{m-1}
\frac{\Gamma(a+j c) \Gamma(b +jc) \Gamma((j+1)c)}
{\Gamma(a +b +(n+j-1)c) \Gamma(c)}
\label{15Veq}
\end{eqnarray}
Here the integration is performed over domain:
$$0 \le t_1 \le t_2 \le \ldots \le t_m \le 1$$
Consider for example $m=2$,  $$V_{\lambda} \otimes V_{2 \Lambda_1}
\longrightarrow V_{\lambda+2 h_2}.$$
$$
a=\frac{(\lambda,-\alpha)}{\kappa}$$
$$b=1+\frac{(2 \Lambda_1, -\alpha_1)}{\kappa} =1-\frac{2}{\kappa}$$
$$c=\frac{1}{2}\frac{ (-\alpha_1,
-\alpha_1)}{\kappa}=\frac{1}{\kappa}$$
Then the Selberg integral is equal:
\begin{eqnarray}
\int \int
(t_1 t_2)^{\frac{(\lambda,
-\alpha)}{\kappa}}(1-t_1)^{-\frac{2}{\kappa}}
(1-t_2)^{-\frac{2}{\kappa}} (t_1-t_2)^{\frac{2}{\kappa}}\frac{
dt_1}{t_1}
\frac{ dt_2}{t_2} \nonumber \\ 
=
\frac{\Gamma(-\frac{(\lambda,\alpha)}{\kappa})
\Gamma(-\frac{(\lambda, \alpha)}{\kappa} +\frac{1}{\kappa})
\Gamma(1-\frac{2}{\kappa})\Gamma(1-\frac{1}{\kappa})
\Gamma(\frac{1}{\kappa})\Gamma(\frac{2}{\kappa})}
{\Gamma(-\frac{(\lambda,\alpha)}{\kappa} +1-\frac{1}{\kappa})
 \Gamma(-\frac{(\lambda,\alpha)+1}{\kappa}) \Gamma(\frac{1}{\kappa})
\Gamma(\frac{1}{\kappa})}
\label{15s1}
\end{eqnarray}
Here the integration is performed over $0 \le t_1 \le t_2 \le 1$.
Passing to another system of contours : contour for $t_1$ starts
and ends at $1$ and encloses $0$ anticlockwise, contour for $t_2$
starts and adds at $1$ encloses $0$ and contour for $t_1$ anticlockwise,
 we get the extra factor:
\begin{eqnarray}
(e^{2 \pi i (-\frac{(\lambda,\alpha)}{\kappa})}-1)
(e^{2 \pi i (-\frac{(\lambda,\alpha)}{\kappa}+\frac{1}{\kappa})}-e^{2
\pi i \frac{1}{\kappa}}
+e^{2 \pi i
(-\frac{(\lambda,\alpha)}{\kappa}+\frac{1}{\kappa})}-1)\nonumber \\
=
e^{\frac{\pi i}{\kappa}} 2 \cos (\frac{\pi}{\kappa}) e^{\pi i \frac{\lambda_2-\lambda_1+1}{\kappa}} (2 i)
\sin(\pi(\frac{\lambda_2-\lambda_1+1}{\kappa}))
 e^{\pi i \frac{(\lambda_2-\lambda_1)}{\kappa}}
(2 i) \sin(\pi (\frac{\lambda_2-\lambda_1}{\kappa})) 
\label{15s2} 
\end{eqnarray}
Multiplying  formulas (\ref{15s1}) by (\ref{15s2}) and using 
$$
\Gamma(x) \Gamma(1-x) = \frac{\pi}{\sin (\pi x)}
$$
we get :
\begin{eqnarray}
\frac{e^{2 \pi i(\frac{\lambda, -\alpha}{\kappa} +\frac{1}{\kappa})}
  \Gamma(1-\frac{1}{\kappa})^2 (2 \pi i)^2}
{\Gamma(\frac{(\lambda,\alpha)}{\kappa}+1-\frac{1}{\kappa})
\Gamma(\frac{(\lambda,\alpha)}{\kappa}+1)
\Gamma(\frac{(\lambda, -\alpha)}{\kappa}+1)
\Gamma(\frac{(\lambda, -\alpha)}{\kappa} +1-\frac{1}{\kappa})}
\label{15s3}
\end{eqnarray}
 The value (\ref{15s3}) coincides with the corresponding value provided by
theorem \ref{15Sel}.
Phase factor $A$ takes into account the fact that first  we had $0 < |t_1| < |z_1| < |t_2| < |z_2| $
and then $0 < |t_1| < |t_2| < z_1 = z_2 =1 $, i.e.
$z_1$ goes through $t_2$ which results in the phase factor
$$A=e^{-\frac{\pi i}  {\kappa}}$$. Also ,
only the leading asymptotic coefficient is essential, while the value at
the unity is equal to $1$ in this example. This phenomena is explained
in the next section.

\section{Particular  case}
The following particular case deserves special attention.
Let $i_1=i_2=i_3=\ldots=i_N=n$.
$$
V_{\lambda} \stackrel{\otimes V_{\Lambda_1}} \longrightarrow
V_{\lambda+h_n} \stackrel{\otimes V_{\Lambda_1}} \longrightarrow
V_{\lambda+2 h_n}\stackrel{\otimes V_{\Lambda_1}}\longrightarrow \ldots \stackrel{\otimes V_{\Lambda_1}}
\longrightarrow
V_{\lambda +N h_n}
$$ 
Then
\begin{eqnarray}
\int \prod t_{ij}^{\frac{(\lambda, -\alpha_j)}{\kappa}}
\prod (z_i -t_{i1})^{\frac{(\Lambda_1, -\alpha_1)}{\kappa}}
\prod (t_{i_1 j_1}-t_{i_2 j_2})^{\frac{(\alpha_{j_1},\alpha_{j_2})}
{\kappa} }\frac{dt_{11}}{t_{11}} \ldots
\frac{dt_{N,n-1}}{t_{N,n-1}}\nonumber \\
=const (z_1 z_2 \ldots z_N)^{\frac{(\lambda,e_n-e_1) +(n-1)}{\kappa}}
\label{13Veq10} 
\end{eqnarray}
Again, the  leading asymptotic coefficient is given by the formula
(\ref{13Veq3}).

 And the Selberg type integral in this case is equal to the leading
asymptotic coefficient (see also the previous section).
Consider $N=n=2$.
Then
\begin{eqnarray}
\int (t_1 t_2)^{\frac{(\lambda, -\alpha)}{\kappa}-1}
(z_1-t_1)^{-\frac{1}{\kappa}}(z_2-t_1)^{-\frac{1}{\kappa}}
(t_2-z_1)^{-\frac{1}{\kappa}}(z_2-t_2)^{-\frac{1}{\kappa}}
(t_2-t_1)^{\frac{2}{\kappa}} dt_1 d t_2\nonumber \\=
const (z_1 z_2)^{\frac{\lambda_2-\lambda_1+1}{\kappa}}
\label{15PR1}
\end{eqnarray}

The constant  is given  by (\ref{13Veq3}).
As a corollary one obtains the following identity for
$\Gamma$-functions for any integer $M$( $M>0$):
\begin{eqnarray}
\sum_{m_1+m_2+m_3=M >0 }
\frac{\Gamma(\frac{(\lambda, -\alpha)}{\kappa} +m_1 +m_3)}
{\Gamma(\frac{(\lambda,-\alpha)}{\kappa} +m_1 +m_3 +1-
\frac{1}{\kappa})} \nonumber \\
\times 
\frac{\Gamma(\frac{(\lambda,\alpha)}{\kappa}+m_2 + m_3)}
{\Gamma(\frac{(\lambda,\alpha)}{\kappa}+m_2+m_3+1 -\frac{1}{\kappa})}
\nonumber \\
\times 
\frac{\Gamma(\frac{1}{\kappa}+m_1) \Gamma(\frac{1}{\kappa}+m_2)
\Gamma(\frac{-2}{\kappa}+m_3)}
{m_1 ! m_2 ! m_3 ! }
= 0
\label{}
\end{eqnarray}

\begin{remark}
Now one can collapse some adjacent variables $z_i$ and obtain an
extension of formula (\ref{13Veq10}).
\end{remark}

{\bf Acknowledgements.}  We acknowledge interesting discussions with 
 I.~Gelfand, V.~Brazhnikov, S.~Lukyanov and A.~Zamolodchikov.

\end{document}